\documentclass[final]{siamltex}

\usepackage{amssymb}
\usepackage{amsmath}
\usepackage{amsfonts}

\DeclareMathOperator \cl {cl\,}

\title{Vector variational principle}

\author{Ewa M. Bednarczuk\thanks{Cardinal Stefan Wyszy\'nski University, 
ul. Dewajtis 5, Warszawa, Poland and Systems Research Institute, 
Polish Academy of Sciences, ul. Newelska 6, 01-447 Warszawa, 
Poland ({\tt Ewa.Bednarczuk@ibspan.waw.pl}).}
        \and Dariusz Zagrodny\thanks{Cardinal Stefan Wyszy\'nski University,  
ul. Dewajtis 5, Warszawa, Poland ({\tt d.zagrodny@uksw.edu.pl}).}}

\begin{document}
\large

\maketitle

\begin{abstract}
 A vector variational principle is proved. 
\end{abstract}



\pagestyle{myheadings}
\thispagestyle{plain}
\markboth{E. M. BEDNARCZUK AND D. ZAGRODNY}{ VECTOR VARIATIONAL PRINCIPLE}

\section{Introduction}
Let  $K\subset Y$ be a closed convex cone  in a locally convex space $Y$ and let
$f:X\longrightarrow Y$ be a mapping defined on  a complete metric space $X$.
Let $D\subset K$ be a  closed convex  bounded subset of $K$
 such that  $0\not \in \cl (D+K)$.
 We show that under mild assumptions on $f$ and $D$ for every $x\in X$ there is $\bar x \in X$ such that 
\begin{equation}\label{nv}
\Big ( f(\bar x)-K\Big)\cap \Big( f(z)+d(z,\bar x)D\Big)=\emptyset \text { for }z\in X\setminus \{\bar x\}.
\end{equation}
If $Y=\mathbb R$, $K=\langle 0,\infty)$, $D=\{\epsilon \}$, $\epsilon >0$, then \eqref {nv} takes the form
$$
f(z)+\epsilon d(z,\bar x)>f(\bar x) \text { whenever }z\not = \bar x.
$$
This is the variational inequality from  Ekeland's Variational Principle (EVP)  \cite {ekeland,ekeland2,phelps1}.
Generalizations of 
EVP to metric spaces are given e.g. in \cite{cammaroto, georgiev, georgiev2, oettli}, 
to locally convex  spaces and to general topological spaces see e.g. \cite{fang}. Thus \eqref{nv} can be regarded as an extension of EVP to vector-valued mappings.

  EVP is a powerful tool with many applications 
in optimization, control theory, subdifferential calculus,  nonlinear analysis, 
global analysis and mathematical economy. Therefore,  several  formulations
of EVP  for vector-valued  and set-valued mappings 
are proved  e.g. in 
\cite{baomardukhovich, bednarczuk, finet1, finet2,gopfert1, gopfert2, isac1, isac2, nemeth, tammer1}. 

However, the common feature of those formulations is their 'directional' character,
more precisely, instead of \eqref{nv} it is proved that
\begin{equation}\label{nv1}
\Big ( f(\bar x)-K\Big)\cap \Big( f(z)+d(z,\bar x)k_{0}\Big)=\emptyset \text { for }z\in X\setminus \{\bar x\},
\end{equation}
 with $0\neq k_{0}$  chosen from the ordering cone $K$.  Under additional assumptions on cone $K$,
N\`emeth (cf. Theorem 6.1 of \cite {nemeth}) proved  \eqref{nv1} with $d(z,\bar x)k_{0}$
replaced by $r(z,\bar{x})$, where $r: X\times X\rightarrow K$ is a mapping such that: 
$(i)$ $r(u,v)=0\Leftrightarrow\ u=v$, $(ii)$  $r(u,v)=r(v,u)$, $(iii)$ $r(u,z)\in r(u,v)+r(v,z)-K$
for any $u,v,z\in X$.

As in the case of scalar-valued mappings, the validity of  EVP for vector- or set-valued mappings is usually verified on the basis of  topological arguments  
and the core of  the proofs is   Cantor's theorem. In contrast to that, in the present paper we prove EVP for vector-valued mappings by combining topological and set-theoretic methods. The main set-theoretic tool is Theorem 3.7 of \cite{gajekzagrodny1} providing sufficient conditions for the existence of maximal elements of countably orderable sets \cite [Definition 2.1]{gajekzagrodny1}. 
The application of set-theoretic methods allows us to prove  \eqref{nv}
which reduces to \eqref{nv1} when $D=\{k_{0}\}$, $0\neq k_{0}\in K$.

The organization of the paper is as follows. In Section 2 we present basic
set-theoretical tools which are used in the sequel. In Section 3 we 
recall semicontinuity concepts for vector-valued functions.
Section 4 contains the main result of the paper, namely Theorem \ref{theorem1} 
together with some examples. In Section 5 the distance $d(x,\bar x)$ is estimated in the case when $x$ is an approximate minimal point and $\bar x$   satisfies \eqref {nv}.

\section{Preliminaries}
The following set-theoretic concepts and facts are used in the sequel.
For any nonempty set $X$ and any relation $s\subset X\times X$ by $x\ s\ y$ we mean that $(x,y)\in s$
and we write $x\ s^{*}\ y$ if and only if
 there is a finite number of elements
 $ x_{1},...,x_{n}\in X$
such that 
$$ x=x_{1},\ x_{1}\ s\ x_{2},...,x_{n-1}\ s\ x_{n}=y
$$
Relation $s^{*}$ is the transitive closure of $s$. If $s$ is transitive, then $s=s^{*}$.

An element $x\in X$ is {\em maximal with respect to a relation $s\subset X\times X$}, we say $x$ is $s$-maximal,
if for every $y\in X$
$$
x\ s\ y\ \ \Rightarrow \ y\ s\ x
$$
When $s$ is a partial order, the above definition  coincides 
with the usual definition of maximality, i.e. $x$ is $s$-maximal if for
every $y\in X$
$$
x\ s\ y\ \ \Rightarrow\ \ x=y,
$$
we refer to \cite{gajekzagrodny,gajekzagrodny1,gajekzagrodny2} 
for more information on maximality with respect to non-transitive relations.

The following definition is essential for the results of Section 3.

\begin{definition}[Definition 2.1 in \cite{gajekzagrodny1}]
\label{definition00}
A set $X$ with a relation $ s\subset X\times X$ is {\em countably orderable}
with respect to $s$ if for every nonempty subset $W\subseteq X$ the existence
of a well ordering relation $\mu$ on $W$ such that
\begin{equation}
\label{orderable}
v\ \mu \ w\ \Rightarrow\ v\ s^{*}\ w\ \ \mbox{for every}\ u,v\in W,\ u\neq w
\end{equation}
implies that  $W$ is at most countable. 
\end{definition}

 Now we are in a position to state  the existence  theorem. 
Its proof has been originally given in \cite {gajekzagrodny} 
(see also Proposition 3.2.8, p. 90 of \cite{gopfert1}).
 
\begin{theorem}
\label{theorem000}
Let $X$ be countably orderable set by a relation $s\subseteq X\times X$.
Assume that for any sequence $(x_{i})\subset X$ satisfying
$$
x_{i}\ s\ x_{i+1}\ \ \mbox{for}\ \ i\in\mathbb{N}
$$
there are a subsequence $(x_{i_{k}})\subset (x_{i})$ and an element $x\in X$ such that
\begin{equation}
\label{th000}
x_{i_{k}}\ s\ x\ \ \mbox{for all } k\in\mathbb{N}.
\end{equation}
Then there exists an $s^{*}$-maximal element of $X$.
\end{theorem}

\begin{corollary}
\label{corollary1}
Let $X$ be countably orderable set by a transitive relation $s\subseteq X\times X$.
Assume that for any sequence $(x_{i})\subset X$ satisfying
$$
x_{i}\ s\ x_{i+1}\ \ \mbox{for}\ \ i\in\mathbb{N}
$$
there are a subsequence $(x_{i_{k}})\subset (x_{i})$ and an element $x\in X$ such that
\begin{equation}
\label{th000-1}
x_{i_{k}}\ s\ x\ \ \mbox{for all } k\in\mathbb{N}.
\end{equation}
Then there exists an $s$-maximal element of $X$.
\end{corollary}

\section{Semicontinuous vector-valued functions}

Let $X$ be a topological  space and let $Y$ be a locally convex space, 
i.e. $Y$ is a linear topological space with a local base consisting of convex neighborhoods 
of the origin, see \cite {holmes}. Let $K\subset Y$ be a closed convex cone in $Y$. For any $x,y\in Y$ define
$$
x\le_{K} y\Leftrightarrow y-x\in K.
$$
We say that a subset $D\subset Y$ is {\em semi-complete} 
if every Cauchy sequence contained in $D$ has a limit in $D$, see \cite {holmes}.
According to \cite{gajekzagrodny} we say that a function $f:X\rightarrow Y$
is {\em monotonically semicontinuous with respect to $K$} (msc)
at $x\in X$ if for every sequence $(x_{i})\subset X$, $x_{i}\rightarrow x$,
satisfying 
$$
f(x_{i+1})\le_{K} f(x_{i})\ \ \ i\in\mathbb{N}
$$
we have $f(x)\le_{K} f(x_{i})$ for $i\in \mathbb{N}$.

We say that $f$ is  msc on $X$ if $f$ is msc at any $x\in X$. To the best of our knowledge the class of msc functions was first introduced 
in \cite[ p.674]{nemeth} via nets. The sequential definition given above is more 
adequate for our purposes. Monotonically semicontinuous functions were   
used in \cite{gajekzagrodny} in order to prove a general form of Weierstrass's theorem for  vector-valued functions and in \cite{nemeth}  to formulate vector  variational principle  (see also Corollary 3.10.19 of \cite{gopfert1}).

We say that $f$ is $K${\em -bounded} if there exists a bounded subset $M$ of $Y$ such that $f(X)\subset K+M$. The topological closure is denoted by $\cl $.

\section{Vector variational principle}

The following theorem provides vector Ekeland's principle.
\vspace{0.3cm}

\begin{theorem}
\label{theorem1}
Let $X$ be a complete metric space and let $Y$ be a locally convex space.
Let $K\subset Y$ be a closed and convex cone in $Y$ and let $D\subset K$ be a closed semi-complete convex and bounded subset of $K$ such that $0\not \in \cl(D+K)$. 

Let $f:X\rightarrow Y$ be msc with respect to $K$ and $K$-bounded. Then for every $x\in X$
there exists $\bar{x}\in X$ such that
\begin{enumerate}
\item[(i)] \begin{tabular}{l}
				$(f(x)-K)\cap(f(\bar{x})+d(x,\bar{x})D))\neq\emptyset$,\end{tabular} 
\item[(ii)] \begin{tabular}{l}
				$(f(\bar{x})-K)\cap(f(z)+d(z,\bar{x})D)=\emptyset$ for every $z\neq\bar{x}$.\end{tabular}
\end{enumerate}
\end{theorem}
\vspace{0.3cm}

\begin{proof} Let $x\in X$ and
$$
A:=\{v\in X\ :\ (f(x)-K)\cap(f(v)+d(x,v)D))\neq\emptyset\}.
$$
Let $r\subset X\times X$ be a relation defined as follows: for any $u,v\in X$
$$
u\ r\ v\ \Leftrightarrow (f(u)-K)\cap(f(v)+d(u,v)D)\neq\emptyset.
$$
By  the convexity of $D$, the relation $r$ is transitive.
The main step of the proof  consists in applying Theorem \ref{theorem000} 
to  show that $A$ has an $r$-maximal element. 
Observe that $A=\{v\in X: x\ r\ v \}$, and hence any $r$-maximal element of $A$ 
is an $r$-maximal element of $X$. 

Since $0\not \in \cl(D+K)$, by separation arguments, there exists $y^{*}\in Y^*$ such that
$$
\langle y^{*},d\rangle+\langle y^{*},k\rangle>\varepsilon>0
$$
for some $\varepsilon>0$ and any $d\in D$, $k\in K$. Hence, $\inf_{d\in D} \langle y^{*},d\rangle>0$ and
$\langle y^{*},k\rangle\ge 0$ for any $k\in K$. Moreover,
\begin{equation}
\label{eq1-theorem1}
u\ r\ v,\ u\neq v\ \ \Rightarrow\  \  \langle y^{*},f(u)\rangle>\langle y^{*},f(v)\rangle.
\end{equation}
Indeed, if $u\ r\ v$, then $f(u)=f(v)+d(u,v)d+k$, where $d\in D$ and $k\in K$.
Consequently, $\langle y^{*},f(u)\rangle=\langle y^{*},f(v)+d(u,v)d+k\rangle >\langle y^{*},f(v)\rangle$.

We start by showing that $A$ is countably orderable with respect 
to $r$. Let $\emptyset\neq W\subseteq A$ be 
any subset of $A$ well
ordered by a relation $\mu$ satisfying (\ref{orderable}). Then for any $u,v\in W$, $u\neq v$
$$
u\ \mu\ v\ \ \Rightarrow\ u\ r\ v\ \Rightarrow\ \langle y^{*},f(u)\rangle>\langle y^{*},f(v)\rangle.
$$
Thus, $y^{*}\circ f(W)\subset \mathbb{R}$ is well ordered by the relation $'>'$
and therefore $f(W)$ is at most countable. This entails that $W$ is countable 
since $y^{*}\circ f$ is a one-to-one mapping on $W$.

Now we  show that \eqref {th000-1} holds for $A$, i.e. for any sequence $(x_{n})\subset A$ 
$$
\forall_{n\in\mathbb{N}}\ x_{n}\ r\ x_{n+1}\ \Rightarrow\ \ \exists_{x_{0}\in A}\ 
\forall_{n\in\mathbb{N}}\ x_{n}\ r\ x_{0}.
$$
By the definition of $r$, for each $n\in\mathbb{N}$
$$
x_{n}\ r\ x_{n+1} \ \ \Leftrightarrow\ \ f(x_{n})-f(x_{n+1})=k_{n}+d(x_{n},x_{n+1})d_{n}\in K,
$$
where $d_{n}\in D$ and $k_{n}\in K$. Moreover, in view of the $K$-boundedness of $f$, for any $m\in\mathbb{N}$
$$
f(x_{1})=\begin{array}[t]{l}
f(x_{1})-f(x_{2})+f(x_{2})-.....-f(x_{m+1})+f(x_{m+1})\\
=f(x_{m+1})+\sum_{i=1}^{m}k_{i}+\sum_{i=1}^{m}d(x_{i},x_{i+1})d_{i}\\
\subset M+K+\sum_{i=1}^{m}d(x_{i},x_{i+1})d_{i}.
\end{array}
$$
Hence, for each $m\in\mathbb{N}$ we have
$$
\langle y^{*},f(x_{1})\rangle\ge\inf_{z\in M}\langle y^{*},z\rangle+
\sum_{i=1}^{m}d(x_{i},x_{i+1})\inf_{d\in D}\langle y^{*},d\rangle.
$$
By the boundedness of $M$ and the fact that $\inf_{d\in D}\langle y^{*},d\rangle>0$, the sequence
$\left(\sum_{i=1}^{m}d(x_{i},x_{i+1})\right)$ is bounded from above, hence
 the series
$$
\sum_{i=1}^{\infty}d(x_{i},x_{i+1})
$$
 converges.  By the boundedness of $D$,   
the sequence $$\sum_{i=1}^{m}d(x_{i},x_{i+1})d_{i}$$ is 
a Cauchy sequence, so by the semi-completness of $D$  the sequence 
$$
\sum_{i=1}^{m}d(x_{i},x_{i+1})d_{i}
$$
converges to a point from $D$ when $m\longrightarrow  \infty $.  Moreover,
$$
\bar{d_{n}}=\sum_{i=n}^{\infty}\frac{d(x_{i},x_{i+1})}{\sum_{i=n}^{\infty}d(x_{i},x_{i+1})}d_{i}
$$ 
is well defined for any $n\in\mathbb{N}$.
In view of  the completeness of $X$, $(x_{n})$ converges to a certain $x_{0}\in X$.
By msc of $f$, $f(x_{n})-f(x_{0})\in K$. 

Now we show that for all $n\in\mathbb{N}$
\begin{equation}
\label{eq-10}
f(x_{n})-f(x_{0})-d(x_{n},x_{0})\bar{d_{n}}\in K.
\end{equation}
Observe that $\bar{d_{n}}\in D$. Because, if it were $\bar{d_{n}}\notin D$, then
in view of the closedness of $D$, it would be $(\bar{d_{n}}+U)\cap  D=\emptyset$
for some neighborhood $U$ of $0$  and for all $m$ sufficiently large, we would get
$$
\sum_{i=n}^{m}\frac{d(x_{i},x_{i+1})}{\sum_{i=n}^{m}d(x_{i},x_{i+1})}d_{i}\in\bar{d_{n}}+U .
$$
The latter, however, would contradict the convexity of $D$. This proves that $\bar{d_{n}}\in D$.

By the definition of $r$, for any $n\in\mathbb{N}$ and $m\ge n$ we obtain 
$$
f(x_{n})=\begin{array}[t]{l}
f(x_{n})-f(x_{n+1})+f(x_{n+1})-.....-f(x_{m+1})+f(x_{m+1})\\
=f(x_{m+1})+\sum_{i=n}^{m}k_{i}+\sum_{i=n}^{m}d(x_{i},x_{i+1})d_{i}.
\end{array}
$$
Since $f(x_{m+1})-f(x_{0})\in K$ this gives
$$
f(x_{n})-f(x_{0})-
d(x_{n},x_{m+1})\left(\frac{\sum_{i=n}^{m}d(x_{i},x_{i+1})d_{i}}{\sum_{i=n}^{m}d(x_{i},x_{i+1})}\right)\in K.
$$
Passing to the limit with $m\rightarrow+\infty$ and taking into account the closedness of $K$ we get (\ref{eq-10}).
 By Theorem \ref{theorem000}, there exists an $r$-maximal element $\bar{x}\in A$. 
Thus, $(i)$ holds for $\bar{x}$. We show that $\bar{x}$ satisfies $(ii)$.
 Since $r$ is transitive,
 $\bar{x}$ is also an $r$-maximal element of $X$. 
If $(ii)$ were not satisfied, we would have $\bar x r z$ for some $z\not =\bar x$ 
and, by the $r$-maximality,  $zr\bar x $. Consequently, by \eqref {eq1-theorem1},  
$\langle y^*,f(\bar x)\rangle > \langle y^*,f(z)\rangle$ and 
$\langle y^*,f(z)\rangle > \langle y^*,f(\bar x)\rangle$, a contradiction.
\end{proof}
\vspace{0.3cm}

Let us note that by repeating the arguments of the proof we can  replace the conclusions $(i)$ and $(ii)$ by
\begin{enumerate}
\item[(i')] \begin{tabular}{l}
				$(f(x)-K)\cap(f(\bar{x})+d(x,\bar{x})(D+K))\neq\emptyset$,\end{tabular} 
\item[(ii')] \begin{tabular}{l}
				$(f(\bar{x})-K)\cap(f(z)+d(z,\bar{x})(D+K))=\emptyset$ for every $z\neq\bar{x}$,\end{tabular}
\end{enumerate}
 where $D$ is as in Theorem \ref{theorem1}. Let us also observe that if $Y$ is a Banach space,  the closedness of $D$ implies its semi-completness, and  $0\not \in \cl (D+K)\ \Leftrightarrow\ d(D+K,0)>0$, 
where $d(D+K,0)$ stands for the distance of $0$ from the set $D+K$. 
In the two examples below it is shown that it is not difficult to verify the inequality $d(D+K,0)>0$.

{\em Examples.
\begin{enumerate}
\item Let $Y$ be a real Banach space, $\varphi$ be from the dual space to $Y$ and $\alpha>0$ be given. If $K$ is contained in a Bishop-Phelps cone generated by $\varphi$ and  $\alpha>0$ i.e. 
\begin{eqnarray*}
K\subset \mathcal{K}_\alpha := \{z\in Y\mid\varphi(z) \ge \alpha\|z\|\}
\end{eqnarray*}
and $d(D,0)>0$, $D\subset K$, then $d(D+K,0)>0$. Indeed, let us observe that for every $d\in D$, $k\in K$ we have 
$$
\| \varphi\| \| d+k\|\ge \varphi(d+k)\ge \varphi(d)\ge \alpha \|d\|\ge \alpha d(D,0)>0.
$$

\item If $K$ is scalarized by a norm according to the idea of Rolewicz \cite{rolewicz} (i.e. $u-v\in K$ implies $\|v\|\le \|u\|$) and $d(D,0)>0$, then $d(D+K,0)>0$.
\end{enumerate}}
\vspace{0.3cm}

\section{Approximate solutions}

For any $d\in D$ put
$\tilde{f}_{d}(z):=f(z)+d(z,\bar{x})d.$
By conclusion $(ii)$ of Theorems \ref{theorem1}, 
$$
(\tilde{f}_{d}(\bar{x})-K)\cap \tilde{f}_{d}(X)= \tilde{f}_{d}(\bar{x})
$$
and $\bar{x}$ is unique in the sense that
$(\tilde{f}_{d}(\bar{x})-K)\cap \tilde{f}_{d}(X\setminus\{\bar{x}\})=\emptyset$.
Thus, for any $d\in D$, $\bar{x}$ is a unique (in the above sense) minimal solution to
problem
$$
\begin{array}{ll}
(P)&K-\min \{\tilde{f}_{d}(x)\ :\ x\in X\}.
\end{array}
$$

Let $D\subset K$ and $0\notin D$, $\varepsilon>0$. We say that $x\in X$ is an 
{\em $\varepsilon-$approximate solution with respect to $D$} for $(P)$ with the objective $f$ if 
$$
(f(x)-\varepsilon D-K)\cap f(X)=\emptyset.
$$
For $K$ being a closed convex pointed cone 
the above definition was given by N\`emeth \cite{nemeth}.

\begin{theorem}
\label{theorem1a}
Let $X$ be a complete metric space and let $Y$ be locally convex space.
Let $K\subset Y$ be a closed and convex cone in $Y$ 
and let $D\subset K$ be a closed semi-complete convex and bounded subset of 
$K$ such that $0\notin \cl (D+K)$. 

Let $f:X\rightarrow Y$ be msc with respect to $K$ and $K$-bounded. Then for every $x\in X$,
$\varepsilon>0$ and $\lambda>0$ there exists $\bar{x}\in X$ such that
\begin{enumerate}
\item[(i)] \begin{tabular}{l}
				$(f(x)-K)\cap(f(\bar{x})+\varepsilon d(x,\bar{x})D))\neq\emptyset$,\end{tabular} 
\item[(ii)] \begin{tabular}{l}
				$(f(\bar{x})-K)\cap(f(z)+\varepsilon d(z,\bar{x})D)=\emptyset$ for every $z\neq\bar{x}$,\end{tabular}
\item[(iii)] Moreover, if $x$ is an $\varepsilon\lambda$-approximate solution with respect to $D$, then
$$
d(x,\bar{x})<\lambda.
$$
\end{enumerate}
\end{theorem}
\vspace{0.3cm}

\begin{proof}  To get the statements $(i)$ and $(ii)$ it is
enough to repeat the proofs of Theorem \ref{theorem1}  with the metric $d(\cdot,\cdot)$
replaced by the metric $\varepsilon d(\cdot,\cdot)$.

Now we  prove $(iii)$.  By $(i)$,
$$
f(x)=f(\bar{x})+\varepsilon d(x,\bar{x})\bar{d}+\bar{k},\ \ \text{where}\ \bar{d}\in D,\ \ \bar{k}\in K.
$$
If it were $d(x,\bar{x})\ge\lambda$, then
$$
f(x)=f(\bar{x})+\varepsilon \frac{d(x,\bar{x})}{\lambda}\lambda\bar{d}+\bar{k}
\in f(\bar{x})+\varepsilon \lambda\bar{d}+ K
$$
which would contradict the fact that $x$ is an $\varepsilon\lambda$-approximate solution.
\end{proof}

\section{Comments} For $D=\{k_{0}\}$ with  $0\neq k_{0}\in K$,  
Theorem \ref{theorem1a} reduces to the results proved e.g. in \cite{baomardukhovich, gopfert1}
and the references therein. Moreover, the proofs given in those references essentially
work for $D=\{k_{0}\}+K$ and only  minor changes are required. 
Let us notice, however, that Theorem \ref{theorem1a} can be applied e.g. to
$D=\{(x_{1},x_{2})\ :\ x_{1}+x_{2}=1,\ x_{1}\ge 0,\ x_{2}\ge 0\}\subset K= \mathbb{R}^{2}_{+}\subset\mathbb{R}^{2}$, 
and $D$ cannot be represented in the form $\{k_{0}\}+K$, where $0\neq k_{0}\in K$.

Generalization of Theorem \ref{theorem1a} in the spirit of $\cite{nemeth}$, 
where $d(z,x)D$ is replaced by a set-valued
mapping $r\colon X\times X\hspace{0.05in}$\parbox{0.18in}{\baselineskip=0.4\baselineskip$
\rightarrow$\\$\rightarrow$}$K$ are  conceivable. The resulting EVP would require some
stringent assumptions on $K$ which we managed to avoid in Theorem \ref{theorem1}
and Theorem \ref{theorem1a}.


\begin{thebibliography}{99}

\bibitem{baomardukhovich}
{\sc T.~Q.~Bao and B.~S.~Mordukhovich}, {\em Variational principles for set-valued
mappings with applications to multiobjective optimization}, 
Control and Cybernetics, 36 (2007), pp.~531--562.

\bibitem{bednarczuk}
{\sc E.~M.~Bednarczuk and ~M.~Przyby\l a}, {\em The vector-valued variational principle in Banach
spaces ordered by cones with nonempty interiors}, 
SIAM J. on Optimization, 18 (2007), pp.~907--913.

\bibitem{cammaroto}
{\sc F.~Cammaroto and A.~Chinni}, {\em A complement to Ekeland's variational principle in Banach spaces}, 
Bull. Pol. Acad. Sci. Math., Volume 44, Number 1, (1996), pp.~29--33.

\bibitem{ekeland}
{\sc I.~Ekeland}, {\em On the variational principle}, 
J. Math. Anal. Appl., 47 (1974), pp.~324--353.

\bibitem{ekeland2}
{\sc I.~Ekeland}{\em Nonconvex minimization problem}, Bull. Amer. Math. Soc. 1 (1979),pp.~443--474.

\bibitem{fang}
{\sc J.~X.~Fang}, {\em The variational principle and fixed points theorems in certain topological spaces}, 
J. Math. Anal. Appl., 208 (1996), pp.~389--412.

\bibitem{finet1}
{\sc C.~Finet}, {\em Variational principles in partially 
ordered Banach spaces}, J. Nonlinear Convex Anal., 
Volume 2, Number 2, (2001), pp.~167--174.

\bibitem{finet2}
{\sc C.~Finet, L.~Quarta and C.~Troestler}, 
{\em Vector-valued variational principles}, 
Nonlinear Anal., 52 (2003), pp.~197--218.

\bibitem{gajekzagrodny}
{\sc L.~Gajek, D. ~Zagrodny}, {\em Weierstrass theorem for monotonically semicontinuous
functions}, Optimization, 29(1994), pp.~199--203

\bibitem{gajekzagrodny2}
{\sc L.~Gajek, D. ~Zagrodny},
{\em Existence of maximal points with respect to ordered bipreference relations},
Journal of Optimization Theory and Applications, 70 (1991), pp.~355--364

\bibitem{gajekzagrodny1}
{\sc L.~Gajek, D. ~Zagrodny},
{\em Countably orderable sets and their applications in optimization},
Optimization, 26 (1992), pp.~287--301

\bibitem{georgiev}
{\sc P.~G.~Georgiev}, {\em Parametric Ekeland's variational principle}, Appl. Math. Lett., 14 (2001), pp.~691-696

\bibitem{georgiev2}
{\sc P.~G.~Georgiev}, {\em The strong Ekeland variational principle, the strong drop theorem and applications}, J. Math. Anal. Appl., 131 (1988), pp.~1-21

\bibitem{gopfert1}
{\sc A.~G\"{o}pfert, H.~Riahi, Ch.~Tammer and C.~Zalinescu}, 
{\em Variational methods in partially ordered spaces}, 
Springer-Verlag, New York, 2003.

\bibitem{gopfert2}
{\sc A.~G\"{o}pfert, Ch.~Tammer and C.~Zalinescu}, 
{\em On the vectorial Ekeland's variational principle 
and minimal point theorems in product spaces}, 
Nonlinear Anal., 39 (2000), pp.~909--922.

\bibitem {holmes}{\sc R. B. Holmes}{\em Geometric functional analysis and its applications}, Springer-Verlag, New York, 1975. 

\bibitem{isac1}
{\sc G.~Isac},
{\em Nuclear cones in product spaces, Pareto efficiency and Ekeland-type variational principle in locally convex spaces},
Optimization, Volume 53, Number 3, (2004), pp.~253-268.

\bibitem{isac2}
{\sc G.~Isac and Ch.~Tammer},
{\em Nuclear and full nuclear cones in product spaces: Pareto efficiency and Ekeland-type variational principle},
Positivity, 9 (2005), pp.~511-539.

\bibitem{nemeth}
{\sc A.B.N\'emeth}, {\em A nonconvex vector minimization problem},
Nonlinear Analysis, Theory and Applications, 10 (1986), pp.~669-678.

\bibitem{oettli}
{\sc W.~Oettli and M.~Th\'{e}ra}, {\em Equivalents of Ekeland's principle},
Bull. Austral. Math. Soc., 48 (1993), pp.~385-392.

\bibitem{phelps1}
{\sc R.~Phelps},
{\em Convex Functions, Monotone Operators and 
Differentiability}, Lecture Notes in Mathematics 1364, 
Springer Verlag, Berlin--Heidelberg, 1989.

\bibitem{phelps2}
{\sc R.~Phelps},
{\em Support cones in Banach spaces and their applications},
 Adv. Math., 13 (1974), pp.~1--19.
 
\bibitem{rolewicz}
{\sc S.~Rolewicz}, {\em On a norm scalarization in infinite dimensional
Banach Spaces}, Control and Cybernetics 4 (1975), pp.~85--89


\bibitem{tammer1}
{\sc Ch.~Tammer},
{\em A generalization of Ekeland's variational principle}, 
Optimization, 25 (1992), pp.~129--141.


\end{thebibliography}
\end{document}